\documentclass[11pt]{amsart}
\usepackage[a4paper,margin=2.5cm]{geometry}
\usepackage[T1]{fontenc}
\usepackage[utf8]{inputenc}
\usepackage{amsmath,amssymb,amsthm,enumitem}
\usepackage{tikz}

\newtheorem{theorem}{Theorem}[section]
\newtheorem{lemma}[theorem]{Lemma}
\newtheorem{proposition}[theorem]{Proposition}
\newtheorem{cor}[theorem]{Corollary}
\theoremstyle{definition}
\newtheorem{definition}[theorem]{Definition}
\newtheorem{example}[theorem]{Example}
\newtheorem{remark}[theorem]{Remark}

\newcommand{\Z}{\mathbb Z}
\newcommand{\N}{\mathbb N}
\newcommand{\F}{\mathbb F}
\newcommand{\C}{\mathbb C}
\DeclareMathOperator{\rk}{rk}\DeclareMathOperator{\im}{Im}

\newcommand{\fonction}[5]{\begin{array}{cccc}
#1: & #2 & \longrightarrow & #3 \\
    & #4 & \longmapsto & #5 \end{array}}

\begin{document}
\title{Three commuting matrices whose lifts never commute}
\author{Charles De Clercq}
\begin{abstract}
We construct pairwise commuting matrices over $\F_p$ which have no pairwise commuting lifts to $\Z/p^2\Z$. Equivalently, there are mod $p$ representations of the absolute Galois group of the field of iterated Laurent series $\C((x_1))\cdots((x_d))$ which do not lift modulo $p^2$.
\end{abstract}
\maketitle

\section{Introduction}

Let $M_1,\dots,M_d$ be pairwise commuting matrices in $M_n(\F_p)$. Can one always find matrices $\tilde M_1,\dots,\tilde M_d\in M_n(\Z/p^2\Z)$, reducing to the $M_l$ modulo $p$, which still commute pairwise?

This very elementary question is related to Galois representations. Let $K=\C((x_1))\cdots((x_d))$ be a field of iterated Laurent series, so that its absolute Galois group is $\widehat{\Z}^d$. A Galois representation $\mathrm{Gal}(\overline K/K)\to\mathrm{GL}_n(\F_p)$ is nothing but $d$ commuting invertible matrices, and it lifts to $\mathrm{GL}_n(\Z/p^2\Z)$ if and only if these matrices have commuting lifts.

The lifting of Galois representations from $\F_p$ to $\Z/p^2\Z$ is a rich subject. Khare and Serre asked whether every continuous representation $\mathrm{Gal}(\overline k/k)\to\mathrm{GL}_n(\F_p)$, for an arbitrary field $k$, lifts modulo $p^2$. The answer is positive when $n=2$ \cite{Khare}, when $p=2$ and $n\le4$ \cite{DeClercqFlorence}, and when $k$ is a local field \cite{Boeckle}. The first counterexamples to statements of this kind were obtained recently by Florence \cite{FlorenceNote} and by Merkurjev and Scavia, who construct in \cite{MerkurjevScavia,MerkurjevScavia2} mod $p$ representations which do not lift modulo $p^2$ in every dimension $n\ge3$ for $p$ odd, and $n\ge5$ for $p=2$. To the best of our knowledge, all previous counterexamples involve fields with non-abelian absolute Galois group.

In this paper we construct new counterexamples, this time for abelian absolute Galois groups, and in a completely elementary way, without any use of Galois cohomology. We construct three commuting matrices\footnote{Two commuting matrices always have commuting lifts (see Section~\ref{sec:2comm}).} over $\F_3$, of size $190$, with no commuting lifts to $\Z/9\Z$, and four commuting matrices over $\F_2$, of size $123$, with no commuting lifts to $\Z/4\Z$. This gives examples of mod $p$ Galois representations of $\C((x_1))\cdots((x_d))$ which do not lift to $\Z/p^2\Z$, for $p=3$ and $d\geq3$, and for $p=2$ and $d\geq4$. We then provide, through cyclic configurations on a cuspidal cubic, three commuting matrices with no commuting lifts for $p\geq 5$ satisfying an arithmetic condition, which holds for all primes lesser than $10^9$ (it remains open for $p=2$).

\section{Preliminaries}

\subsection{Notation}
All rings are commutative with unit. In the sequel, we fix a prime $p$, a finite field $k$ of characteristic $p$ and the ring $W_2(k)$ of Witt vectors of length $2$ (for $k=\F_p$, $W_2(k)$ is simply $\Z/p^2\Z$). The ring $W_2(k)$ is local with maximal ideal $pW_2(k)$, and multiplication by $p$ induces an isomorphism $k\xrightarrow{\sim}pW_2(k)$. We often write $\overline{\cdot}$ for the reduction $W_2(k)\to k$ of elements in $W_2(k)$-modules.

We fix an integer $d\ge1$ and always endow $\Z^d$ with the \emph{product order}: $u\le v$ if and only if $u_l\le v_l$ for $l=1,\dots,d$. We write $e_1,\dots,e_d$ for the standard basis of $\Z^d$ and, for $s\le t$, we denote by $[s,t]=\{w\in\Z^d: s\le w\le t\}$ the box with corners $s$ and $t$.

\subsection{The lifting problem}

\begin{definition}\label{def:lift}
Let $V$ be a $k$-vector space of finite dimension $n$ and $u_1,\dots,u_d$ pairwise commuting endomorphisms of $V$. A \emph{lift} of $(V;u_1,\dots,u_d)$ to $W_2(k)$ is a free $W_2(k)$-module $V_2$ with an isomorphism $\alpha\colon V_2\otimes_{W_2(k)}k\to V$, together with pairwise commuting endomorphisms $\tilde u_1,\dots,\tilde u_d$ of $V_2$ such that $\alpha\circ(\tilde u_l\otimes k)=u_l\circ\alpha$ for $l=1,\dots,d$.

We say that $(V;u_1,\dots,u_d)$ \emph{lifts} if it admits such a lift.
\end{definition}

Note that $(V;u_1,\dots,u_d)$ lifts if and only if in some basis of $V$ (equivalently in any basis, by Nakayama's lemma), the matrices $N_1,\dots,N_d$ of $u_1,\dots,u_d$ lift to pairwise commuting matrices $\tilde N_1,\dots,\tilde N_d\in M_n(W_2(k))$.

\section{Two commuting matrices always lift}\label{sec:2comm}

In this section only, $k$ is an arbitrary field and $R$ is a ring together with a surjective ring morphism $\pi\colon R\to k$. We write $a\mapsto\bar a$ for $\pi$, as well as for the induced surjections $R[x]\to k[x]$ and $M_n(R)\to M_n(k)$. Two commuting matrices in $M_n(k)$ always lift to commuting matrices in $M_n(R)$ (\cite[Lemma~6]{JP12} for uniserial rings of length two, and \cite[Lemma~2.8]{PSS15}, due to Singla, for any ring with a surjective morphism onto $k$). We recall a short proof for completeness, and deduce that all mod $p$ Galois representations of $\C((x))((y))$ lift to $\Z_p$.

\begin{lemma}\label{lem:2liftdiv}
    Let $P$, $Q$ in $R[x]$ be monic, with $P|Q$ or $Q|P$. A morphism of $k[x]$-modules $\varphi\colon k[x]/(\bar{Q})\to k[x]/(\bar{P})$ is the reduction of a morphism of $R[x]$-modules $\Phi\colon R[x]/(Q)\to R[x]/(P)$.
\end{lemma}

\begin{proof}
A morphism of $R[x]$-modules $\Phi\colon R[x]/(Q)\to R[x]/(P)$ is determined by $\Phi(1)$ which can be any element such that $Q\cdot\Phi(1)=0$ (and similarly over $k$).

If $P\mid Q$, then any lift of $\varphi(1)$ in $R[x]/(P)$ for $\Phi(1)$ gives rise to a lift $\Phi$ of $\varphi$, since $Q$ (being a multiple of $P$) kills this lift. 

If $Q\mid P$, write $P=QE$ with $E$ monic. Since $\bar Q\cdot\varphi(1)=0$ in $k[x]/(\bar P)$ and $k[x]$ is an integral domain, $\varphi(1)$ is the class of $\bar E S$ for some $S\in k[x]$. Lift $S$ to $S'\in R[x]$ and take $\Phi(1)$ to be the class of $ES'$. It defines a lift $\Phi$ since it is killed by $Q$ (as $QES'=PS'$) and reduces to $\varphi(1)$.
\end{proof}

\begin{theorem}\label{thm:lift2endo}
    For any finite dimensional $k$-vector space $V$ and commuting endomorphisms $u$ and $v$ of $V$, $(V;u,v)$ lifts to $R$ through $\pi$ (Definition \ref{def:lift} replacing $W_2(k)$ by $R$).
\end{theorem}

\begin{proof}
View $V$ as a $k[x]$-module, $x$ acting by $u$. As $u$ and $v$ commute, $v$ defines a morphism of $k[x]$-modules. The structure theorem of finitely generated $k[x]$-modules yields non-constant monic polynomials $d_1\mid d_2 \mid \cdots \mid d_r$ in $k[x]$ and an isomorphism of $k[x]$-modules $$V\simeq \bigoplus_{i=1}^rk[x]/(d_i),$$
so that $v$ is given by a family of $k[x]$-linear maps $v_{ij}:k[x]/(d_j)\to k[x]/(d_i)$, with $1\leq i,j \leq r$.

Choose inductively monic lifts $D_1\mid D_2 \mid\cdots \mid D_r$ of $d_1,\dots,d_r$ by picking for $D_1$ a lift of $d_1$, and then set $D_i=D_{i-1}Q_i$, where $Q_i$ is a monic lift of the quotient $d_i/d_{i-1}$. Setting
$$\tilde V:=\bigoplus_{i=1}^rR[x]/(D_i),$$
we obviously get a lift of $V$, on which multiplication $\tilde{u}$ by $x$ lifts $u$. Moreover since either $D_i\mid D_j$ or $D_j \mid D_i$, Lemma \ref{lem:2liftdiv} states that each of the $v_{ij}$ lifts to morphisms of $R[x]$-modules $\tilde{v}_{ij}:R[x]/(D_j)\to R[x]/(D_i)$. These define an $R[x]$-linear endomorphism $\tilde{v}$ of $\tilde V$ which reduces to $v$ and commutes with $\tilde{u}$ (the multiplication by $x$). 
\end{proof}

As explained right after Definition \ref{def:lift}, Theorem \ref{thm:lift2endo} implies that two commuting matrices in $M_n(k)$ lift to commuting matrices in $M_n(R)$, if $R$ is a local ring with residue field $k$. The general case can be obtained with the following trick.

\begin{cor}\label{cor:lift2mat}
Let $k$ be a field, $R$ an arbitrary ring with a surjective ring morphism $\pi\colon R\to k$. Two commuting matrices in $M_n(k)$ lift to two commuting matrices in $M_n(R)$.
\end{cor}

\begin{proof}
Let $A$ and $B$ in $M_n(k)$ be two commuting matrices, and denote by $u$ and $v$ the associated endomorphisms of $V=k^n$. By (the proof of) Theorem \ref{thm:lift2endo}, $u$ and $v$ lift over $\tilde V=\bigoplus_{i=1}^rR[x]/(D_i)$ to commuting endomorphisms $\tilde{u}$ and $\tilde{v}$. Denote by $\tilde C$ and $\tilde E$ the matrices of $\tilde u$ and $\tilde v$ in the basis $\tilde{\mathcal{B}}$ of $\tilde V$ given by powers of $x$, on each block. These matrices commute by definition and their images $C$ and $E$ through $\pi$ are the matrices of $u$ and $v$ in the basis $\mathcal{B}$ of $V\simeq \bigoplus_{i=1}^rk[x]/(d_i)$ given by monomials on blocks as well.

Denote by $P$ the change-of-basis matrix from the canonical basis of $V$ to $\mathcal{B}$, so that $A=PCP^{-1}$ and $B=PEP^{-1}$. Lifting each coefficient of $P$ through $\pi$, choose a matrix $\tilde{P}\in M_n(R)$ reducing to $P$, and pick $\delta\in R$ such that $\pi(\delta)=\det(P)^{-1}$. Set $$\tilde A=\delta\cdot \tilde P \cdot\tilde C\cdot \mathrm{Adj}(\tilde P)\mbox{  and  }\tilde B=\delta\cdot\tilde P\cdot \tilde E\cdot \mathrm{Adj}(\tilde P)$$
where $\mathrm{Adj}(\cdot)$ denotes the adjugate matrix (transpose of the cofactor matrix). As $\tilde C$ and $\tilde E$ commute and since $\mathrm{Adj}(\tilde P)\tilde P=\det(\tilde P)I_n$,  an easy computation shows that $\tilde A \tilde B= \tilde B \tilde A$.
Note that the coefficients of $\mathrm{Adj}(\tilde P)$ are polynomials in the coefficients of $\tilde P$, so that the reduction to $k$ of $\mathrm{Adj}(\tilde P)$ is $\mathrm{Adj}(P)=\det(P)P^{-1}$. It follows that the reduction of $\tilde A$ is $PCP^{-1}=A$ and the reduction of $\tilde B$ is $B$.\end{proof}

\begin{cor}\label{cor:lift2galois}
Let $K=\C((x))((y))$, so that $\mathrm{Gal}(\overline K/K)\simeq\widehat{\Z}^2$. Every Galois representation $\rho\colon\mathrm{Gal}(\overline K/K)\to\mathrm{GL}_n(\F_p)$ lifts to a continuous representation $\tilde\rho\colon\mathrm{Gal}(\overline K/K)\to\mathrm{GL}_n(\Z_p)$, and in particular to $\mathrm{GL}_n(\Z/p^m\Z)$ for every $m\geq1$.
\end{cor}

\begin{proof}
Let $\sigma_1,\sigma_2\in\mathrm{Gal}(\overline K/K)$ correspond to the standard basis of $\widehat\Z^2$. By Corollary~\ref{cor:lift2mat} applied to $R=\Z_p$, the commuting matrices $\rho(\sigma_1),\rho(\sigma_2)$ have commuting lifts $\tilde A_1,\tilde A_2\in M_n(\Z_p)$, invertible since $\Z_p$ is local. For every $m\ge1$, the images of $\tilde A_1,\tilde A_2$ in the finite group $\mathrm{GL}_n(\Z/p^m\Z)$ commute, hence define a continuous representation $\tilde\rho_m$ of $\widehat\Z^2$. As these representations are compatible with the reductions $\mathrm{GL}_n(\Z/p^{m+1}\Z)\to\mathrm{GL}_n(\Z/p^m\Z)$, they give rise to a continuous representation $\tilde\rho:\mathrm{Gal}(\overline K/K)\to\mathrm{GL}_n(\Z_p)$ which lifts $\rho$, since its reduction coincides with $\rho$ on the dense subset $\Z^2$ of $\widehat\Z^2$.
\end{proof}

\section{$\Z^d$-graded modules and placed matrices}

From now on, $k$ is a finite field of characteristic $p$ and we only consider lifts to the local ring $W_2(k)$. We now introduce the basics on $\Z^d$-graded modules and their lifts. These notions are classical, and part of a much richer theory (see for instance \cite{MS} and \cite{Mi}). We shall only need a few definitions and elementary facts.

\begin{definition}
Let $A$ be a ring. A \emph{$\Z^d$-graded module} over $A$ is a family of $A$-modules $(M_w)_{w\in\Z^d}$ together with $A$-linear maps $x_l\colon M_w\to M_{w+e_l}$ ($l=1,\dots,d$) such that $x_lx_m=x_mx_l$ for all $l,m$. For $c=(c_1,\dots,c_d)\in\N^d$ we write $x^c=x_1^{c_1}\circ\cdots \circ x_d^{c_d}\colon M_w\to M_{w+c}$.
\end{definition}

Equivalently, $\bigoplus_{w\in \Z^d} M_w$ is a module over $A[x_1,\dots,x_d]$, graded with respect to the usual multidegree (this is the usual notion of a $\Z^d$-graded module \cite{MS}). It is also the same as a module over the poset $(\Z^d,\le)$ in the sense of \cite[Def.~2.1]{Mi}, with $\varphi_{w,w+c}=x^c$.

\begin{definition}
A $\Z^d$-graded module $M$ over $k$ has \emph{finite length} if $\dim_k M_w < \infty$ for all $w \in \Z^d$ and $M_w=0$ for all but finitely many $w$.
\end{definition}

If $M$ is a $\Z^d$-graded module of finite length over $k$, $V=\bigoplus_{w\in \Z^d}M_w$ is a finite dimensional vector space and $x_1,\dots,x_d$ are commuting nilpotent endomorphisms of $V$. We say that $M$ lifts (to $W_2(k)$) if $(V;x_1,\dots,x_d)$ lifts in the sense of Definition~\ref{def:lift}.

\begin{definition}
Let $M$ be a $\Z^d$-graded module of finite length over $k$. A \emph{$\Z^d$-graded lift} of $M$ is a $\Z^d$-graded module $\tilde M$ over $W_2(k)$, with maps $\tilde x_l\colon\tilde M_w\to\tilde M_{w+e_l}$, such that each $\tilde M_w$ is a free $W_2(k)$-module, together with for any $w\in \Z^d$ an isomorphism $\alpha_w\colon\tilde M_w\otimes_{W_2(k)}k\xrightarrow{\sim}M_w$ such that $\alpha_{w+e_l}\circ(\tilde x_l\otimes k)=x_l\circ\alpha_w$ for all $w$ and $l$ (that is, on each component, $\tilde x_l$ reduces to $x_l$ via $\alpha_w$).
\end{definition}

\begin{lemma}[reduction to $\Z^d$-graded lifts]\label{lem:graded}
A $\Z^d$-graded module $M$ of finite length and dimension $n$ over $k$ lifts to $W_2(k)$ if and only if it has a $\Z^d$-graded lift.
\end{lemma}
\begin{proof}
The data of a $\Z^d$-graded lift $\tilde M$ gives a lift of $M$, with $V_2=\bigoplus_w\tilde M_w$ and commuting endomorphisms $\tilde{x}_1,\dots,\tilde{x}_d$. Conversely, choose a basis $\mathcal{B}_w$ of each $M_w$ and concatenate them into a basis $\mathcal{B}$ of $V$. In this basis, for any $i=1,\dots, d$, the endomorphism $x_i$ of $V$ has a block matrix $N_i$. Pick a lift $\tilde{N}_i$ of $N_i$ in $M_n(W_2(k))$, with the same block configuration.

Multiplication by $p$ induces a morphism of $W_2(k)$-modules $\iota\colon M_n(k)\to M_n(W_2(k))$, $\iota(Z)=p\hat Z$ for any lift $\hat Z$ of $Z$ (well defined since $p^2=0$, injective with image $pM_n(W_2(k))$). Since $[N_i,N_j]$ is trivial, the commutator $[\tilde N_i,\tilde N_j]$ reduces to $0$, so that $[\tilde N_i,\tilde N_j]=\iota(C_{ij})$ for a unique $C_{ij}\in M_n(k)$. This yields for any $X_1,\dots,X_d$ in $M_n(k)$ the following direct computation, using $p^2=0$: 
$$[\tilde N_i+\iota(X_i),\tilde N_j+\iota(X_j)]=[\tilde N_i,\tilde N_j]+\iota\bigl([N_i,X_j]-[N_j,X_i]\bigr)=\iota\bigl(C_{ij}+[N_i,X_j]-[N_j,X_i]\bigr).$$
So, $\iota$ being injective, the matrices $\tilde N_i+\iota(X_i)$ and $\tilde N_j+\iota(X_j)$ with $i\neq j$ commute if and only if
$$(\ast)\qquad [N_i,X_j]-[N_j,X_i]=-C_{ij},\mbox{ for all $i\neq j$}.$$

By assumption, there are pairwise commuting matrices $L_1,\dots,L_d$ in $M_n(W_2(k))$ reducing to $N_1,\dots, N_d$. We have $L_i=\tilde N_i + \iota(X_i)$, where the $X_i$ satisfy $(\ast)$.

Decompose matrices of $M_n(k)$ into blocks $(w',w)$ according to $V=\bigoplus_{w\in \Z^d}M_w$. Given an element $Z\in M_n(k)$ and $v\in\Z^d$, denote by $Z^{(v)}\in M_n(k)$ the matrix obtained from $Z$ by keeping only the blocks $(w+v,w)$. Note that for all $Z_1$, $Z_2$ in $M_n(k)$ and $v_1$, $v_2$ in $\Z^d$, $Z_1^{(v_1)}Z_2^{(v_2)}$ vanishes out of the blocks $(w+v_1+v_2,w)$, $w\in \Z^d$. Since $N_i$ only has blocks $(w+e_i,w)$, we thus have $[N_i,Z]^{(v+e_i)}=[N_i,Z^{(v)}]$ for every $Z$ and every $v$. Moreover $C_{ij}$ only has blocks $(w+e_i+e_j,w)$, because $\tilde N_i$ and $\tilde N_j$ have the same blocks as $N_i$ and $N_j$ (hence are homogeneous of respective degrees $e_i$ and $e_j$). Keeping only the blocks $(w+e_i+e_j,w)$ in $(\ast)$ thus gives
$$[N_i,X_j^{(e_j)}]-[N_j,X_i^{(e_i)}]=-C_{ij}.$$
Therefore the $X_i^{(e_i)}$ satisfy $(\ast)$ and the matrices $\tilde L_i:=\tilde N_i+\iota(X_i^{(e_i)})$ pairwise commute. They reduce to $N_i$ and only have blocks $(w+e_i,w)$, so that setting $\tilde M_w=W_2(k)^{\mathcal B_w}$, the morphisms $\tilde x_i\colon \tilde M_w\to \tilde M_{w+e_i}$ given by the blocks $(w+e_i,w)$ of $\tilde L_i$ define a $\Z^d$-graded lift (with $\alpha_w$ given by the basis $\mathcal B_w$).
\end{proof}

\begin{remark}
Lemma~\ref{lem:graded} is specific to $W_2(k)$: there the lifting problem is linear, which
allows one to keep only the homogeneous components of the right degree. Over $W_m(k)$ for
$m\ge3$, or over ramified rings, quadratic terms appear, and a graded module might lift
without lifting as a graded module.
\end{remark}

\begin{remark}\label{rem:schmiermann}
By work of Schmiermann \cite{Sch23}, configurations of points and lines in $\mathbb P^2$ which do not lift to $\Z/p^2\Z$ give, via Lemma~\ref{lem:graded}, commuting matrices without commuting lifts: for $d=3$ when the incidence graph of the configuration is a string graph \cite[Cor.~26]{Sch23}, and always for $d\ge4$ \cite[Thm.~30]{Sch23}. Theorem~\ref{thm:p2} is an explicit such example, of dimension $123$, coming from the Fano plane.
\end{remark}

\subsection{Placed matrices}

We now build $\Z^d$-graded modules of finite length from small matrices whose rows and columns are placed at points of $\Z^d$. This construction is not new: it is a special case of Miller’s fringe and flange presentations of modules over posets, encoded by monomial matrices \cite[\S\,3.2, \S\,5.4]{Mi}, which go back to the monomial matrices for free and injective modules of \cite[\S\,1.4, \S\,11.3]{MS}. We give an elementary description, enough for our purpose.

\begin{definition}\label{def:placed}
A \emph{placed matrix} over a ring $A$ consists of a matrix $G=(g_{ij})\in M_{a\times b}(A)$ and points $t_1,\dots,t_a\in\Z^d$ (the \emph{sites of the rows}) and $s_1,\dots,s_b\in\Z^d$ (the \emph{sites of the columns}), such that $g_{ij}=0$ unless $s_j\le t_i$. For $w\in\Z^d$ put
$$I(w)=\{i\in\{1,\dots,a\}, w\leq t_i\},\qquad J(w)=\{j\in\{1,\dots,b\}, s_j\leq w \},$$
and call the submatrix $G_w=(g_{ij})_{i\in I(w),\,j\in J(w)}$ the \emph{block at $w$}.
\end{definition}

\begin{remark} If $I(w)$ or $J(w)$ is empty, $G_w$ is an empty matrix, of rank $0$.
\end{remark}  

\begin{definition}\label{def:module}
To a placed matrix $G$ over $A$ we associate the $\Z^d$-graded module $M(G;A)$ defined as follows: $M(G;A)_w\subseteq A^{I(w)}$ is the submodule spanned by the columns of $G_w$, and $x_l\colon M(G;A)_w\to M(G;A)_{w+e_l}$ is the restriction of the natural projection $\pi_l:A^{I(w)}\to A^{I(w+e_l)}$.
\end{definition}

In other words, $M(G;A)_w$ is spanned by the columns of $G_w$. For $1\le l\le d$, the block $G_{w+e_l}$ has fewer rows and more columns than $G_w$, and $x_l$ removes from each column of $G_w$ its entries from rows $i\in I(w)\setminus I(w+e_l)$, which gives the corresponding column of $G_{w+e_l}$. These maps commute, since $x_lx_m$ and $x_mx_l$ are both given by the projection $A^{I(w)}\to A^{I(w+e_l+e_m)}$.

\begin{remark}\label{rem:every}
Every $\Z^d$-graded module of finite length over $k$ is of the form $M(G)$ for some placed matrix $G$, whose columns correspond to homogeneous generators and whose rows correspond to a basis of the dual of the socle. This is the finite length case of Miller's fringe presentations \cite{Mi} (the construction is the one in the proof of Theorem~\ref{thm:criterion}, where it is carried out over $W_2(k)$).
\end{remark}
 
\begin{proposition}\label{prop:tuple}
Let $G$ be a placed matrix over $k$. Then $M(G)=M(G;k)$ has finite length, its nonzero components lie in the union of the boxes $[s_j,t_i]$, and its dimension is
\[\dim(M(G))=\sum_{w\in\Z^d}\rk(G_w).\]
In particular, with $n=\dim(M(G))$, the matrices $N_1,\dots,N_d\in M_n(k)$ of $x_1,\dots,x_d$ in a basis of $\bigoplus_wM(G)_w$ commute and are nilpotent.
\end{proposition}
\begin{proof}
By definition $\dim(M(G)_w)=\rk(G_w)$. If $M(G)_w\ne0$, there is some $i\in I(w)$ and $j\in J(w)$, i.e.\ $s_j\le w\le t_i$. The $x_l$ commute, raise the degree, and only finitely many components are nonzero.
\end{proof}
 
\begin{example}
\begin{enumerate}
    \item If $a=b=d=1$, $G=(1)$, $s_1=s$, $t_1=t$, $M(G)_w=k$ for $w\in[s,t]$ and $0$ otherwise. The matrix $N_1$ of $x_1$ on $M(G)$ is thus a nilpotent Jordan block, and every $\Z$-graded module of finite length is a direct sum of such modules. This is the bar code of persistence \cite[Ex.~3.19--3.20]{Mi}.

    \item Let $d=2$, and let $G$ be the placed matrix over $k$
\[
G=\begin{pmatrix}1&1\\1&2\end{pmatrix},\qquad s_1=(0,1),\ s_2=(1,0),\qquad t_1=(3,2),\ t_2=(2,3).
\]
The figure shows $r_w=\rk(G_w)=\dim M(G)_w$ at each point $w$ of the box $[(0,0),(3,3)]$. Outside of this box, all components are trivial.
\begin{center}
\begin{tikzpicture}[scale=1.1]
\foreach \x in {0,...,3} \foreach \y in {0,...,3} \fill[gray!40] (\x,\y) circle (1.5pt);
\def\r{{{0,1,1,1},{1,2,2,1},{1,2,2,1},{1,1,1,0}}}
\foreach \x in {0,...,3} \foreach \y in {0,...,3} {
  \pgfmathtruncatemacro{\v}{\r[\x][\y]}
  \ifnum\v>0 \node[fill=white,inner sep=1pt] at (\x,\y) {$\v$};\fi}
\draw[thick,blue] (0,1) circle (7pt); \node[blue,left=6pt] at (0,1) {$s_1$};
\draw[thick,blue] (1,0) circle (7pt); \node[blue,below=6pt] at (1,0) {$s_2$};
\draw[thick,red] (2.75,1.75) rectangle (3.25,2.25); \node[red,right=6pt] at (3,2) {$t_1$};
\draw[thick,red] (1.75,2.75) rectangle (2.25,3.25); \node[red,above=6pt] at (2,3) {$t_2$};
\draw[->] (4.2,-0.6) -- (5,-0.6) node[right] {$x_1$};
\draw[->] (4.2,-0.6) -- (4.2,0.2) node[above] {$x_2$};
\end{tikzpicture}
\end{center}
For instance, for $w=(1,1)$, $I(w)=J(w)=\{1,2\}$ so that $G_w=G$ and $r_w=2$, while for $w=(3,1)$, $G_w=(1\ 1)$ has rank $1$. The morphism $x_1^2\colon M(G)_{(1,1)}\to M(G)_{(3,1)}$ is the projection on the first coordinate and $\dim M(G)=18$.
\end{enumerate} 
\end{example}

\subsection{A lemma on minors}
\begin{lemma}\label{lem:minors}
Let $B\in M_{m,q}(W_2(k))$, let $\bar B$ be its reduction and set $r=\rk(\bar B)$. The submodule of $W_2(k)^m$ spanned by the columns of $B$ is a free direct summand of rank $r$ if and only if all minors of $B$ of order $r+1$ vanish.
\end{lemma}
\begin{proof}
Since $\bar B$ has rank $r$, some minor of order $r$ of $\bar B$ is nonzero, so the corresponding minor of $B$ is a unit in $W_2(k)$. Gaussian elimination with unit pivots, after placing this minor at the top-left corner, gives $P\in\mathrm{GL}_m(W_2(k))$ and $Q\in\mathrm{GL}_q(W_2(k))$ such that
$$B'=PBQ=\begin{pmatrix}I_r&0\\0&D\end{pmatrix} \quad\text{with $\bar D=0$ since $\bar B$ has rank $r$}.$$
We have $\im (B')=P(\im (B))$ and $B$, $B'$ have the same ideal of minors of order $r+1$, we may thus assume that $B=B'$. Now $\im (B)=W_2(k)^r\oplus\im D\subseteq W_2(k)^r\oplus W_2(k)^{m-r}$, and it remains to observe that both conditions of the lemma are equivalent to $D=0$ (for the condition on minors, this is clear). If $D=0$, then $\im (B)=W_2(k)^r$ is obviously a free direct summand of rank $r$. Conversely, if $\im (B)$ is free of rank $r$, then $\im (B)\otimes_{W_2(k)}k$ has dimension $r$. Since $p\,\im (D)=0$, this tensor product is $k^r\oplus\im (D)$, hence $\im (D)=0$.
\end{proof}

\section{Lifting modules of placed matrices}\label{sec:criterion}

The following result shows that the module $M(G)$ of a placed matrix lifts to $W_2(k)$ if and only if $G$ itself has a lift, as a placed matrix with the same sites, and which preserves the rank of the blocks $G_w$, for all $w\in\Z^d$.
 
\begin{theorem}\label{thm:criterion}
Let $G$ be a placed matrix over $k$, with sites $t_1,\dots,t_a$ and $s_1,\dots,s_b$, and let $r_w=\rk(G_w)$ for $w\in\Z^d$. The $\Z^d$-graded module $M(G)$ lifts to $W_2(k)$ if and only if there is a placed matrix $\tilde G$ over $W_2(k)$, with the same sites, which reduces to $G$ and such that, for every $w\in\Z^d$, all minors of order $r_w+1$ of $\tilde G_w$ vanish.
\end{theorem}
 
\begin{proof}
First, assume there is a placed matrix $\tilde G$ over $W_2(k)$ sharing the sites, lifting $G$ and satisfying the condition on minors. For any $w\in \Z^d$, the $w$-component of the module $M(\tilde G;W_2(k))$ is spanned by the columns of $\tilde{G}_w$ and is a free direct summand of rank $r_w$ of $W_2(k)^{I(w)}$ by Lemma \ref{lem:minors}. In particular the natural morphism $M(\tilde G;W_2(k))_w\otimes_{W_2(k)}k\to W_2(k)^{I(w)}\otimes_{W_2(k)}k=k^{I(w)}$ is injective, and its image is spanned by the reductions of the columns of $\tilde G_w$, hence is simply $M(G)_w$. This identification is compatible with the morphisms $x_1,\dots,x_d$, which are projections on both sides. The module $M(\tilde G;W_2(k))$ is thus a $\Z^d$-graded lift of $M(G)$.

If $M(G)$ lifts, by Lemma \ref{lem:graded} it has a $\Z^d$-graded lift $(\tilde M;\tilde x_1,\dots ,\tilde x_d;(\alpha_w)_{w\in \Z^d})$. We now show that this lift can be identified with the module of a placed matrix over $W_2(k)$ satisfying the condition on minors.

\emph{Construction of $\tilde G$.} For $1\leq i \leq a$, let $\lambda_i:M(G)_{t_i}\longrightarrow k$
be the $i$-th coordinate of $M(G)_{t_i}\subset k^{I(t_i)}$. As $\tilde M_{t_i}$ is free, we can lift $\lambda_i$ into a morphism of $W_2(k)$-modules $\tilde\lambda_i:\tilde M_{t_i}\longrightarrow W_2(k)$.

Now for $1\leq j \leq b$, lift the $j$-th column of $G_{s_j}$ (an element of $M(G)_{s_j}$) to an element $\tilde \gamma_j\in \tilde M_{s_j}$. Setting $\tilde{g}_{ij}=\tilde\lambda_i(\tilde{x}^{t_i-s_j}\tilde{\gamma}_j)$ if $s_j\leq t_i$ and $\tilde{g}_{ij}=0$ otherwise, we get a placed matrix $\tilde{G}$ over $W_2(k)$ with same sites as $G$. The reduction of $\tilde{g}_{ij}$ to $k$ is the $i$-th coordinate of the $j$-th column of $G_{t_i}$, that is $g_{ij}$, so that $\tilde G$ reduces to $G$.

\emph{Condition on minors.} Fix $w\in \Z^d$ and consider the morphism of $W_2(k)$-modules
$$\fonction{\varphi_w}{\tilde{M}_w}{W_2(k)^{I(w)}}{m}{(\tilde\lambda_i(\tilde{x}^{t_i-w}m))_{i\in I(w)}}.$$
Reduced to $k$, the $i$-th component of $\varphi_w(m)$ is the $i$-th coordinate of the projection of $\bar m\in M(G)_w$ to $M(G)_{t_i}$, that is, the $i$-th coordinate of $\bar m$. Hence $\varphi_w$ reduces to the inclusion $M(G)_w\subseteq k^{I(w)}$. As $\tilde M_w$ is free of rank $r_w$ over the local ring $W_2(k)$ and the reduction of $\varphi_w$ is injective, $\varphi_w$ is a split injection whose image is a free direct summand of $W_2(k)^{I(w)}$, of rank $r_w$.

We now show that for any $w\in \Z^d$, $\varphi_w(\tilde{M}_w)=M(\tilde G;W_2(k))_w$. Inclusion $\supseteq$ is clear as $\varphi_w(\tilde x^{\,w-s_j}\tilde\gamma_j)=\bigl(\tilde\lambda_i(\tilde x^{\,t_i-s_j}\tilde\gamma_j)\bigr)_{i\in I(w)}$ is the $j$-th column of $\tilde G_w$, for all $j\in J(w)$. For the other inclusion, by Nakayama's lemma the $\tilde x^{\,w-s_j}\tilde\gamma_j$, with $j\in J(w)$, span $\tilde{M}_w$, since their reductions to $k$ are the columns of $G_w$. The $\varphi_w(\tilde x^{\,w-s_j}\tilde\gamma_j)$, $j\in J(w)$, are the columns of $\tilde{G}_w$ and span $\varphi_w(\tilde{M}_w)$. We have shown that for any $w\in \Z^d$, the $W_2(k)$-module $M(\tilde{G};W_2(k))_w$ generated by the columns of $\tilde{G}_w$ is a free direct summand of $W_2(k)^{I(w)}$, of rank $r_w$. As the reduction of $\tilde{G}_w$ is of rank $r_w$, Lemma \ref{lem:minors} implies that all minors of order $r_w+1$ of $\tilde{G}_w$ vanish.
\end{proof}
 
\section{Four commuting matrices for $p=2$}\label{sec:p2}

\begin{theorem}\label{thm:p2}
Let $G=(g_{ij})\in M_{3,4}(\F_2)$ be the placed matrix
\[
G=\begin{pmatrix}0&1&1&1\\1&0&1&1\\1&1&0&1\end{pmatrix},\qquad
\begin{array}{ll}
t_1=(7,7,7,5), & s_1=(5,6,5,4),\\
t_2=(8,7,6,5), & s_2=(6,4,6,4),\\
t_3=(8,7,7,4), & s_3=(6,7,4,3),\\
 & s_4=(7,5,4,4),
\end{array}
\]
Then $M(G)$ has dimension $123$, and it does not lift to $W_2(\F_2)=\Z/4\Z$. In other words, the matrices $N_1,\dots,N_4\in M_{123}(\F_2)$ of $x_1,\dots,x_4$ on $M(G)$ commute pairwise, and there are no pairwise commuting matrices in $M_{123}(\Z/4\Z)$ reducing to them.
\end{theorem}

Since $s_j\le t_i$ for all $i,j$, $G$ is a placed matrix. Denoting by $M(G)=M(G;\F_2)$ the associated module, a direct computation gives $\dim(M(G))=\sum_{w\in\Z^4}\rk(G_w)=123$.

The idea is the following: we write down conditions on the minors of a $3\times4$ matrix which can be satisfied over $\F_2$ but not over $\Z/4\Z$ (Lemma~\ref{lem:fano}). The matrix $G$ satisfies them over $\F_2$, and its sites are chosen so that each of these minors comes from a block $G_w$. By Theorem~\ref{thm:criterion}, a lift of $M(G)$ would give a matrix over $\Z/4\Z$ satisfying the same conditions, which is impossible.

\begin{lemma}\label{lem:fano}
Let $A$ be a commutative ring, and let $H=(h_{ij})\in M_{3,4}(A)$ be a matrix such that
\begin{enumerate}
\item $h_{11}=h_{22}=h_{33}=0$, and the nine other entries are units of $A$;
\item the three $2\times2$ submatrices with rows $\{2,3\}$ and columns $\{1,4\}$, rows $\{1,3\}$ and columns $\{2,4\}$, rows $\{1,2\}$ and columns $\{3,4\}$ have determinant $0$;
\item the $3\times3$ submatrix with columns $\{1,2,3\}$ has determinant $0$.
\end{enumerate}
Then $2=0$ in $A$. In particular there is no such matrix over $\Z/4\Z$.
\end{lemma}

\begin{proof}
Conditions $(1)$, $(2)$ and $(3)$ are stable under multiplication of a row or a column by a unit. Dividing the rows $1,2,3$ by $h_{14},h_{24},h_{34}$, then the columns $1,2,3$ by the new entries $h_{21},h_{12},h_{13}$, we may assume that for some units $a$, $b$ and $c$,
$$H=\begin{pmatrix}0&1&1&1\\1&0&c&1\\a&b&0&1\end{pmatrix}.$$
The three conditions (2) give $a=b=c=1$, so the determinant of condition (3) is
$$\begin{vmatrix}0&1&1\\1&0&1\\1&1&0\end{vmatrix}=2,$$
and condition (3) gives $2=0$ in $A$.
\end{proof}

\begin{remark}
Over a field, the lemma says that such matrices only exist in characteristic $2$. This is the classical fact that the Fano matroid is representable only in characteristic $2$ \cite[Prop.~6.4.8]{Oxley}. The matrix $G$ is the so-called fundamental-circuit incidence matrix of the Fano matroid $F_7$ with respect to the basis $\{1,2,3\}$, see \cite[Fig.~6.12(b)]{Oxley}.
\end{remark}

\begin{proof}[Proof of Theorem~\ref{thm:p2}]
For the following seven points $w$, we compute the sets $I(w)$ and $J(w)$, and the rank $r_w=\rk(G_w)$ of each associated matrix over $\F_2$:
$$\begin{array}{llll}
w & I(w) & J(w) & r_w\\ \hline
(5,6,7,5) & \{1\} & \{1\} & 0\\
(8,4,6,5) & \{2\} & \{2\} & 0\\
(8,7,7,3) & \{3\} & \{3\} & 0\\
(8,6,5,4) & \{2,3\} & \{1,4\} & 1\\
(7,5,7,4) & \{1,3\} & \{2,4\} & 1\\
(7,7,4,5) & \{1,2\} & \{3,4\} & 1\\
(6,7,6,4) & \{1,2,3\} & \{1,2,3\} & 2
\end{array}$$
These are exactly the seven submatrices of Lemma~\ref{lem:fano}, with the prescribed ranks.

Assuming that $M(G)$ lifts, Theorem~\ref{thm:criterion} states that there is a placed matrix $\tilde G=(\tilde g_{ij})$ in $M_{3,4}(\Z/4\Z)$, with the same sites and reducing to $G$, such that for every $w\in\Z^4$ all minors of order $r_w+1$ of $\tilde G_w$ vanish. For the seven points of the table, this means: $\tilde g_{11}=\tilde g_{22}=\tilde g_{33}=0$, the three $2\times2$ determinants of condition (2) vanish, and the $3\times3$ determinant of condition (3) vanishes. Moreover the nine other entries of $\tilde G$ reduce to nonzero elements of $\F_2$, hence are units, so that $\tilde G$ would satisfy the conditions of Lemma~\ref{lem:fano} over $\Z/4\Z$, a contradiction.
\end{proof}

\begin{cor}\label{cor:galois2}
Let $K=\C((x_1))((x_2))((x_3))((x_4))$, so that $\mathrm{Gal}(\overline K/K)\simeq\widehat{\Z}^4$. There is a Galois representation $\rho\colon\mathrm{Gal}(\overline K/K)\to\mathrm{GL}_{123}(\F_2)$ which does not lift to $\mathrm{GL}_{123}(\Z/4\Z)$.
\end{cor}

\begin{proof}
Let $\sigma_1,\dots,\sigma_4\in\mathrm{Gal}(\overline K/K)$ correspond to the standard basis of $\widehat{\Z}^4$, and let $N_1,\dots,N_4$ be the matrices of Theorem~\ref{thm:p2}. Since the matrices $1+N_l$ commute pairwise and $\mathrm{GL}_{123}(\F_2)$ is finite, $\rho(\sigma_l):=1+N_l$ defines a continuous representation. If $\tilde\rho$ were a lift of $\rho$ to $\mathrm{GL}_{123}(\Z/4\Z)$, the matrices $\tilde\rho(\sigma_l)-1$ would commute pairwise and reduce to $N_l$, contradicting Theorem~\ref{thm:p2}.
\end{proof}

\begin{remark}[Computer checks]
We could not find an example of a $\Z^3$-graded module of finite length over $\F_2$ which does not lift to $\Z/4\Z$, contrary to the case of odd primes. Computer searches suggest that there may be none.
\end{remark}

\section{Three commuting matrices for $p=3$}\label{sec:p3}

\begin{theorem}\label{thm:p3}
Let $G=(g_{ij})\in M_{3,5}(\F_3)$ be the placed matrix
\[
G=\begin{pmatrix}1&0&2&1&2\\1&1&2&2&0\\0&2&2&2&1\end{pmatrix},\qquad
\begin{array}{ll}
t_1=(3,6,3), & s_1=(2,0,2),\\
t_2=(4,4,4), & s_2=(2,3,-1),\\
t_3=(5,5,2), & s_3=(3,1,0),\\
 & s_4=(1,2,1),\\
 & s_5=(0,4,0),
\end{array}
\]
Then $M(G)$ has dimension $190$, and it does not lift to $W_2(\F_3)=\Z/9\Z$. In other words, the matrices $N_1,N_2,N_3\in M_{190}(\F_3)$ of $x_1,x_2,x_3$ on $M(G)$ commute pairwise, but there are no pairwise commuting matrices in $M_{190}(\Z/9\Z)$ reducing to them.
\end{theorem}

Since $s_j\leq t_i$ for all $i,j$, $G$ is a placed matrix, and denoting by $M(G)=M(G;\F_3)$ the associated module, a direct computation gives $\dim(M(G))=\sum_{w\in \Z^3}\rk(G_w)=190$. The idea of the proof is the same as in Section~\ref{sec:p2}, with Lemma~\ref{lem:fano} replaced by the following lemma.

\begin{lemma}\label{lem:mk}
Let $A$ be a commutative ring, and let $H=(h_{ij})\in M_{3,5}(A)$ be a matrix such that
\begin{enumerate}
\item $h_{12}=h_{25}=h_{31}=0$, and the twelve other entries are units of $A$;
\item the three $2\times2$ submatrices with rows $\{1,2\}$ and columns $\{1,3\}$, rows $\{1,3\}$ and columns $\{4,5\}$, rows $\{2,3\}$ and columns $\{3,4\}$ have determinant $0$;
\item the two $3\times3$ submatrices given by columns $\{1,2,4\}$ and $\{2,3,5\}$ have determinant $0$.
\end{enumerate}
Then there is $u\in A$ satisfying $u^2-u+1=0$. In particular there is no such matrix over $\Z/9\Z$.
\end{lemma}

\begin{proof}
Conditions (1), (2) and (3) are stable under multiplication of a row or a column of $H$ by a unit. Dividing the columns $1,3,4,5$ by $h_{11},h_{13},h_{14},h_{15}$, then the rows $2,3$ by the new entries $h_{21},h_{33}$, then the column $2$ by the new entry $h_{22}$, we may assume that for some units $a,b,c,u,e$,
\[
H=\begin{pmatrix}1&0&1&1&1\\1&1&a&b&0\\0&c&1&u&e\end{pmatrix}.
\]
The conditions (2) give $a=1$, $e=u$ and $b=u$. The second determinant of condition (3), given by the columns $\{2,3,5\}$, is thus
\[
\begin{vmatrix}0&1&1\\1&1&0\\c&1&u\end{vmatrix}=1-c-u,
\]
so that $c=1-u$, and the first determinant given by columns $\{1,2,4\}$ is then
\[
\begin{vmatrix}1&0&1\\1&1&u\\0&1-u&u\end{vmatrix}=u+(1-u)-u(1-u)=u^2-u+1.
\]
Hence $u^2-u+1=0$, an equation which no element of $\Z/9\Z$ satisfies.
\end{proof}

\begin{remark}
Over $\F_3$, the columns of $[I_3\,|\,G]$ form the Möbius--Kantor configuration of eight points and eight lines, the lines corresponding to the eight conditions of Lemma~\ref{lem:mk}. As a matroid, it is $AG(2,3)\backslash e$, one of the ``interesting matroids'' of \cite[Appendix]{Oxley}.
\end{remark}

\begin{proof}[Proof of Theorem~\ref{thm:p3}]
For the following eight points $w$, we compute the sets $I(w)$ and $J(w)$, and the rank $r_w=\rk(G_w)$ of each associated matrix over $\F_3$:
$$\begin{array}{llll}
w & I(w) & J(w) & r_w\\ \hline
(2,6,-1) & \{1\} & \{2\} & 0\\
(0,4,4) & \{2\} & \{5\} & 0\\
(5,0,2) & \{3\} & \{1\} & 0\\
(3,1,3) & \{1,2\} & \{1,3\} & 1\\
(1,5,1) & \{1,3\} & \{4,5\} & 1\\
(4,2,1) & \{2,3\} & \{3,4\} & 1\\
(2,3,2) & \{1,2,3\} & \{1,2,4\} & 2\\
(3,4,0) & \{1,2,3\} & \{2,3,5\} & 2
\end{array}
$$
These are exactly the eight submatrices of Lemma~\ref{lem:mk}, with the prescribed ranks.

Assuming that $M(G)$ lifts, Theorem~\ref{thm:criterion} states that there is a placed matrix $\tilde G=(\tilde g_{ij})$ in $M_{3,5}(\Z/9\Z)$, with the same sites and reducing to $G$, such that for every $w\in\Z^3$ all minors of order $r_w+1$ of $\tilde G_w$ vanish. Applied to the eight points of the table, and since the twelve other entries of $\tilde G$ are units, $\tilde G$ would satisfy the conditions of Lemma~\ref{lem:mk} over $\Z/9\Z$, a contradiction.
\end{proof}

\begin{remark}[Computer checks]
The size of $G$ in Theorem \ref{thm:p3} is minimal: the module of every placed matrix over $\F_3$ with at most $3$ rows and $4$ columns, or at most $4$ rows and $3$ columns, lifts to $\Z/9\Z$, for all choices of sites. Note that the dimension of $M(G)$ is not determined by the size of $G$: we constructed placed matrices in $M_{4,4}(\F_3)$ and $M_{8,8}(\F_3)$ whose modules do not lift to $\Z/9\Z$ either, of respective dimensions $172$ and $146$.
\end{remark}

\begin{cor}\label{cor:galois}
Let $K=\C((x_1))((x_2))((x_3))$, so that $\mathrm{Gal}(\overline K/K)\simeq\widehat{\Z}^3$. There is a Galois representation $\rho\colon\mathrm{Gal}(\overline K/K)\to\mathrm{GL}_{190}(\F_3)$ which does not lift to $\mathrm{GL}_{190}(\Z/9\Z)$.
\end{cor}

\begin{proof}
Same as the proof of Corollary~\ref{cor:galois2}, using Theorem~\ref{thm:p3}.
\end{proof}

\section{Larger primes}\label{sec:larger}

The method used for $p=3$ can be extended to larger primes (through computer search we constructed $\Z^3$-graded modules over $\F_p$ of finite length which do not lift to $\Z/p^2\Z$ for $p=5$, $7$, $11$ and $13$, of respective dimensions $151$, $251$, $221$ and $1372$). We provide a uniform construction for $p\geq 5$, based on a cyclic configuration of points on a cuspidal cubic curve, under an arithmetic condition. We first explain it for $p=11$ (where the placed matrix is not too big).

Throughout, $\zeta \in \overline{\F}_p$ is a root of $x^3+x+1$ of multiplicative order $n$ at least $4$, $k=\F_p(\zeta)$, and $\omega\in W_2(k)$ is the Teichmüller lift of $\zeta$.

\noindent\textbf{A cyclic configuration on a cuspidal cubic.} In the sequel, we denote by $(C_1,\dots,C_m)$ the matrix with columns $C_1,\dots,C_m$. The triples $P(a):={}^t\!\begin{pmatrix}1&a&a^3\end{pmatrix}$ in $k^3$ are the smooth points of the cuspidal cubic $X^2Z=Y^3$ in $\mathbb{P}^2$. Vandermonde computation gives
$$\det(P(y_1),P(y_2),P(y_3))=(y_2-y_1)(y_3-y_1)(y_3-y_2)(y_1+y_2+y_3)$$
so that three distinct such points are collinear if and only if $x_1+x_2+x_3=0$.

Since $\zeta$ has order $n$, the point $P(\zeta^l)$ only depends on $l\in\Z/n\Z$, and these $n$ points are distinct. As $\zeta^l+\zeta^{l+1}+\zeta^{l+3}=\zeta^l(1+\zeta+\zeta^3)=0$, the points $P(\zeta^l)$, $P(\zeta^{l+1})$, $P(\zeta^{l+3})$ are collinear for every $l\in\Z/n\Z$. In other words, the matrix $H_{\zeta}=(P(\zeta^l))_{l\in\Z/n\Z}\in M_{3,n}(k)$ satisfies $\det(P(\zeta^l),P(\zeta^{l+1}),P(\zeta^{l+3}))=0$ for every $l$. The set $\mathrm{Aff}(\zeta)$ defined below is the set of the lifts of $H_{\zeta}$ to $W_2(k)$ which satisfy the same conditions.

\begin{definition}
Let $\mathrm{Aff}(\zeta)$ be the set of matrices $\tilde H=(\tilde h_l)_{l\in\Z/n\Z}$ in
$M_{3,n}(W_2(k))$ which reduce to $H_{\zeta}$ and satisfy $\det(\tilde h_l,\tilde h_{l+1},\tilde h_{l+3})=0$ for every $l\in\Z/n\Z$.
\end{definition}

The set of lifts of $H_{\zeta}$ in $M_{3,n}(W_2(k))$ is an affine $k$-space with direction $M_{3,n}(k)$, since $pW_2(k)\simeq k$. Furthermore, by multilinearity of the determinant and since $p^2=0$ in $W_2(k)$, the conditions $\det(\tilde h_l,\tilde h_{l+1},\tilde h_{l+3})=0$ are affine, so that $\mathrm{Aff}(\zeta)$ itself is an affine $k$-space (possibly empty). The elements of $\mathrm{Aff}(\zeta)$ are the lifts of the configuration, and the next lemma shows that such lifts exist if and only if the Teichmüller lift $\omega$ of $\zeta$ is again a root of $x^3+x+1$.

\begin{lemma}\label{lem:aff}
    $\mathrm{Aff}(\zeta)$ is non-empty if and only if $\omega^3+\omega+1=0$.
\end{lemma}

\begin{proof}
First assume $\mathrm{Aff}(\zeta)$ is non-empty and denote by $M_{\omega}$ the diagonal matrix with entries $1,\omega,\omega^3$. The affine map
$$\fonction{\sigma}{\mathrm{Aff}(\zeta)}{\mathrm{Aff}(\zeta)}{(\tilde h_0,\dots,\tilde h_{n-1})}{(M_{\omega}\tilde h_{n-1},M_{\omega}\tilde h_0,\dots,M_{\omega}\tilde h_{n-2})}$$
given by multiplying columns by $M_{\omega}$ and applying a cyclic permutation is well-defined, since the reduction of $M_{\omega}$ is the diagonal matrix with entries $1,\zeta,\zeta^3$ and since the appropriate minors still vanish. As $p\nmid n$ ($n$ divides $|k|-1$), a point $s\in \mathrm{Aff}(\zeta)$ gives rise to a $\sigma$-invariant point, considering the barycenter $\frac{1}{n}\sum_{j=0}^{n-1}\sigma^j(s)$. We may thus assume that $s=(\tilde h_0,\dots,\tilde h_{n-1})$ is $\sigma$-invariant, hence that $\tilde h_l=M_{\omega}^l\tilde h_0$, for all $l$. As $\tilde h_0$ reduces to $P(1)$, $\tilde h_0={}^t(a,b,c)$ for some invertible elements in $W_2(k)$, hence $\tilde h_l={}^t (a,\omega^l b,\omega^{3l} c)$. Vanishing of the minor of $s$ associated with columns $\tilde h_0,\tilde h_1,\tilde h_3$ gives
$$0=abc\begin{vmatrix}
   1&1&1\\1&\omega&\omega^3\\1&\omega^3&\omega^9 
\end{vmatrix}=abc(\omega-1)(\omega^3-1)(\omega^3-\omega)(\omega^3+\omega+1)$$
and $\omega^3+\omega+1=0$, since $\zeta$ has order at least $4$.

On the other hand if $\omega^3+\omega+1=0$, setting $H_{\omega}=\begin{pmatrix}P(\omega^l)\end{pmatrix}_{l\in \Z/n\Z}$, we get an obvious lift of $H_{\zeta}$ which belongs to $\mathrm{Aff}(\zeta)$ since for any $l$, 
$$\det(P(\omega^l),P(\omega^{l+1}),P(\omega^{l+3}))=\omega^{l}(\omega^{l+1}-\omega^{l})(\omega^{l+3}-\omega^{l})(\omega^{l+3}-\omega^{l+1})(\omega^3+\omega+1)=0.$$
\end{proof}

\noindent\textbf{The case $p=11$.} In $\F_{11}$, $\zeta=2$ satisfies our assumptions, as a root of $x^3+x+1$ of order $10$, whose Teichmüller lift $\omega=2^{11}\in \Z/121\Z$ does not satisfy $\omega^3+\omega+1=0$. Hence $\mathrm{Aff}(2)$ is empty.

\begin{theorem}\label{thm:p11}
Let $G$ be the placed matrix over $\F_{11}$
$$G=\left(\begin{array}{*{13}{c}}
   1&1&1&1&1&1&1&1&1&1&1&1&1\\
   1&1&2&4&4&8&5&5&10&9&7&3&6\\
   1&1&8&9&9&6&4&4&10&3&2&5&7
   \end{array}\right)$$
with $$
\begin{array}{lllll}
s_1=(2,2,12), & s_2=(0,11,5), & s_3=(5,10,2), & s_4=(1,12,1), & s_5=(11,8,0), \\
s_6=(10,5,4), & s_7=(12,0,3), & s_8=(9,1,9), & s_9=(8,6,6), & s_{10}=(7,4,10), \\
s_{11}=(6,7,7), & s_{12}=(4,3,11), & s_{13}=(3,9,8)&\\
\end{array}$$
and $t_1=t_2=t_3=(12,12,12)$ (the coordinate-wise maximum of the $s_j$'s).

Then $M(G)$ has dimension $1640$ and does not lift to $\Z/121\Z$.
\end{theorem}

\begin{remark}
The matrix $G$ above is obtained from $H_{\zeta}$ with $\zeta=2$ by duplicating its first, third and fifth column, that is $P(\zeta^0)$, $P(\zeta^2)$ and $P(\zeta^4)$. For other primes, the same construction works, with sites given by explicit formulas in $n$ (see Theorem~\ref{thm:larger}).
\end{remark}
 
\begin{proof}
Since $s_j$'s are lesser than every $t_i$, $I(w)=\{1,2,3\}$ for every $w$ lesser than $(12,12,12)$. For the following points, we get the respective $J(w)$ and $r_w=\rk(G_w)$:
\[
\begin{array}{lll}
w & J(w) & r_w\\\hline
(2,11,12) & \{1,2\} & 1\\
(11,12,1) & \{4,5\} & 1\\
(12,1,9) & \{7,8\} & 1\\
(10,11,5) & \{2,3,6\} & 2\\
(12,12,3) & \{3,4,5,7\} & 2\\
(11,8,6) & \{5,6,9\} & 2\\
(12,5,10) & \{6,7,8,10\} & 2\\
(9,7,9) & \{8,9,11\} & 2\\
(8,6,11) & \{9,10,12\} & 2\\
(7,9,10) & \{10,11,13\} & 2\\
(6,7,12) & \{1,11,12\} & 2\\
(5,10,11) & \{3,12,13\} & 2\\
(3,12,12) & \{1,2,4,13\} & 2
\end{array}
\]
Assume $M(G)$ lifts, and let $\tilde G$ be given by Theorem~\ref{thm:criterion}. By the first line of the table, columns $1$ and $2$ of the lift $\tilde G$ must be proportional, with a unit factor since they reduce to $P(1)$. The same holds by the second and third line of the table for columns $4,5$ and $7,8$ of $\tilde G$, which reduce to $P(4)$ and $P(5)$ respectively. The matrix $\tilde H=(\tilde h_l)_{l\in\Z/10\Z}$ obtained from $\tilde G$ by removing the first, fourth and seventh column is a lift of the matrix $H_2$ associated to $\zeta=2$. 

Using these proportionalities, the ten last lines of the table, whose blocks have rank $2$, give $\det(\tilde h_l,\tilde h_{l+1},\tilde h_{l+3})=0$ for every $l\in\Z/10\Z$ (for instance, the line $J(w)=\{3,4,5,7\}$ gives $l=1$, since column $7$ is proportional to column $8$). Hence $\tilde H\in\mathrm{Aff}(2)$, which contradicts Lemma~\ref{lem:aff}. The dimension is a long but direct computation.
\end{proof}
 
\noindent\textbf{Generalization to $p\ge5$.} The case $p=11$ above uses the root $\zeta=2$ of order $n=10$. In general, the placed matrix is built from $H_\zeta$ in the same way: its entries depend on $\zeta$, but its size and its sites only depend on the order $n$ of $\zeta$, which we take even, hence at least $10$ by Lemma~\ref{lem:even} below. The only delicate point is the choice of the sites. We found them by a computer search, and the proof of Theorem \ref{thm:larger} shows that they have the required properties for every even $n\ge10$.
 
\begin{lemma}\label{lem:even}
For $p\ge5$, the polynomial $x^3+x+1$ has a root of even order in $\overline\F_p$, and every root of even order has order at least $10$.
\end{lemma}
 
\begin{proof}
The product of the three roots is $-1$, of order $2$, so they cannot all have odd order. If $\zeta$ is a root of even order $n\le8$, then $\zeta^{n/2}=-1$. One easily checks that if $\zeta=-1$, $\zeta^2=-1$ or $\zeta^3=-1$, $\zeta$ is not a root of $x^3+x+1$. In the same way if $\zeta^4=-1$, $\zeta^4+\zeta^2+\zeta=0$ gives $-1=\zeta^4=(1-\zeta)^2=2-3\zeta$ so we would have $\zeta=1$ ($p\neq 3$) which is not a root.
\end{proof}
 
\begin{theorem}\label{thm:larger}
Let $p\ge5$, and assume that $\zeta$ has even order $n$ and that $\omega^3+\omega+1\neq0$ in $W_2(k)$. Then there are three commuting matrices over $\F_p$ which have no commuting lifts over $\Z/p^2\Z$.
\end{theorem}
 
The condition $\omega^3+\omega+1\neq0$ for some root $\zeta$ of even order holds for every prime $5 \leq p < 10^9$. The matrix of Theorem~\ref{thm:p11} corresponds to the case $n=10$. Let $G\in M_{3,n+3}(k)$ be obtained from $H_\zeta$ by duplicating the columns corresponding to the points $P(\zeta^0)$, $P(\zeta^2)$ and $P(\zeta^4)$:
$$G=\bigl(P(\zeta^0),P(\zeta^0),P(\zeta),P(\zeta^2),P(\zeta^2),P(\zeta^3),P(\zeta^4),P(\zeta^4),P(\zeta^5),\dots,P(\zeta^{n-1})\bigr).$$

We denote by $v_1\vee\dots\vee v_r$ the coordinatewise maximum of $v_1,\dots,v_r$ in $\Z^3$.

Now for the sites: set $t_1=t_2=t_3=s_1\vee\dots\vee s_{n+3}$, where
\[
\begin{array}{ll}
s_1=(-2n,-n,3n), & s_2=(-2n-4,n+4,3),\\
s_3=(-2n+5,n+2,-2), & s_4=(-2n-2,n+6,-18),\\
s_5=(-4,n-3,-20), & s_6=(-6,2,2),\\
s_7=(0,-n-12,0), & s_8=(-8,-n-10,2n+4),
\end{array}
\]
and, for $5\le l\le n-1$,
\[
s_{l+4}=\begin{cases}(-2l,\,l-1,\,l-1) & \text{if $l$ is odd,}\\ (-2l,\,-l,\,l+2n) & \text{if $l$ is even.}\end{cases}
\]
\begin{proof}[Proof of Theorem~\ref{thm:larger}]
By Lemma~\ref{lem:even}, $n \geq 10$. We first show that $M(G)$ does not lift to $W_2(k)$, the argument being, mutatis mutandis, the one of Theorem~\ref{thm:p11} (bearing in mind here $\zeta$ may not belong to $\F_p$). The following $w$ are enough (note that $I(w)=\{1,2,3\}$ for all of them):
{\small
\[
\begin{array}{llll}
w & J(w) & r_w & \\\hline
s_1\vee s_2=(-2n,n+4,3n) & \{1,2\} & 1 & \\
s_4\vee s_5=(-4,n+6,-18) & \{4,5\} & 1 & \\
s_7\vee s_8=(0,-n-10,2n+4) & \{7,8\} & 1 & \\
s_2\vee s_3\vee s_6=(-6,n+4,3) & \{2,3,6\} & 2 & l=0\\
s_3\vee s_5\vee s_7=(0,n+2,0) & \{3,5,7\} & 2 & l=1\\
s_5\vee s_6\vee s_9=(-4,n-3,4) & \{5,6,9\} & 2 & l=2\\
s_6\vee s_8\vee s_{10}=(-6,2,2n+6) & \{6,8,10\} & 2 & l=3\\
s_8\vee s_9\vee s_{11}=(-8,6,2n+4) & \{8,9,11\} & 2 & l=4\\
s_{l+4}\vee s_{l+5}\vee s_{l+7}=(-2l,\,l-1,\,l+2n+3) & \{l+4,l+5,l+7\} & 2 & l \text{ odd},\ 5\le l\le n-5\\
s_{l+4}\vee s_{l+5}\vee s_{l+7}=(-2l,\,l+2,\,l+2n) & \{l+4,l+5,l+7\} & 2 & l \text{ even},\ 6\le l\le n-4\\
s_{n+1}\vee s_{n+2}\vee s_1=(-2n+6,n-4,3n) & \{1,n+1,n+2\} & 2 & l=n-3\\
s_{n+2}\vee s_{n+3}\vee s_3=(-2n+5,n+2,3n-2) & \{3,n+2,n+3\} & 2 & l=n-2\\
s_{n+3}\vee s_1\vee s_4=(-2n+2,n+6,3n) & \{1,2,4,n+3\} & 2 & l=n-1
\end{array}
\]}
As for $p=11$, the first three lines provide collinearity of columns $1-2$, $4-5$ and $7-8$ of a presumed lift $M(\tilde G)$ of $M(G)$ to $W_2(k)$. As removing three columns of $G$ yields $H_{\zeta}$, the lift $\tilde G$ would give an element of $\mathrm{Aff}(\zeta)$, contradicting Lemma~\ref{lem:aff}.

It remains to descend to $\F_p$. Let $e=[k:\F_p]$ and $N=\dim_k M(G)$. As an $\F_p$-vector space, $M(G)$ has dimension $eN$, and the multiplication $m_{\zeta}$ by $\zeta$ and $x_1,x_2,x_3$ are commuting $\F_p$-linear endomorphisms. Let $A_1,A_2,A_3\in M_{eN}(\F_p)$ be the matrices of $m_{\zeta}+x_1$, $x_2$, $x_3$ in some $\F_p$-basis, and set $q:=|k|$. Assume they have commuting lifts $\tilde A_1,\tilde A_2,\tilde A_3$ over $\Z/p^2\Z$. Since $m_\zeta^q=m_\zeta$, $x_1$ is nilpotent and $m_\zeta$, $x_1$ commute, we have $A_1^{q^r}=m_\zeta+x_1^{q^r}=m_\zeta$ for $r$ large enough, hence $\tilde A_1^{q^{r+1}}= \tilde A_1^{q^r}\mathrm{~mod}~p$. Since $X=Y\mathrm{~mod}~p$ implies $X^p=Y^p\mathrm{ mod}~p^2$ for commuting $X,Y$, raising to the power $q$ gives $\tilde A_1^{q^{r+2}}=\tilde A_1^{q^{r+1}}$. Thus $\tilde m_{\zeta}=\tilde A_1^{q^{r+1}}$ lifts $m_\zeta$, satisfies $\tilde m_{\zeta}^q=\tilde m_{\zeta}$ and being a polynomial in $\tilde A_1$, it commutes with $\tilde A_2$ and $\tilde A_3$.

Write $W_2(k)=(\Z/p^2\Z)[X]/(\tilde\mu)$, where $\tilde\mu$ is the monic lift of the minimal polynomial $\mu$ of $\zeta$ dividing $X^q-X$, and write $X^q-X=\tilde\mu\tilde g$, with $\tilde g$ reducing to $g$. Since $X^q-X$ is separable, we have $g(\zeta)\neq0$. Hence $\tilde g(\tilde m_\zeta)$, which reduces to $g(m_\zeta)=m_{g(\zeta)}$, is invertible, and $\tilde m_\zeta^q=\tilde m_\zeta$ gives $\tilde\mu(\tilde m_\zeta)=0$. So $X\mapsto\tilde m_\zeta$ makes $(\Z/p^2\Z)^{eN}$ a $W_2(k)$-module. It is free of rank $N$ by Nakayama's lemma, since it is free over $\Z/p^2\Z$ and its reduction modulo $p$ is $M(G)$. As $\tilde A_2$ and $\tilde A_3$ commute with $\tilde A_1$, they commute with $\tilde m_\zeta$, so $\tilde A_1-\tilde m_\zeta$, $\tilde A_2$, $\tilde A_3$ are commuting $W_2(k)$-linear endomorphisms lifting $x_1,x_2,x_3$. This contradicts the first part of the proof.
\end{proof}
 
\begin{remark}
The matrices of Theorem~\ref{thm:larger} have size $eN$, where $e=[k:\F_p]$ and $N=\dim M(G)$. One can show that $N \leq 3(n+3)^3$, but this dimension is generally large. For instance for $p=5$, we have $n=62$ and get $N=257831$, hence matrices of size $eN=773493$ over $\F_5$ (compare with the other example above, of dimension $151$).
\end{remark}

\end{document}